\magnification=1200
\overfullrule=0pt
\centerline {\bf A note on spherical maxima sharing the same Lagrange multiplier}\par
\bigskip
\bigskip
\centerline {BIAGIO RICCERI}\par
\bigskip
\bigskip
{\bf Abstract.}
 In this paper, we establish a general result on spherical maxima sharing the same
Lagrange multiplier of which the following is a particular consequence: Let $X$ be a real Hilbert
space. For each $r>0$, let $S_r=\{x\in X : \|x\|^2=r\}$. Let $J:X\to {\bf R}$ be a sequentially
weakly upper semicontinuous functional which is G\^ateaux differentiable in $X\setminus \{0\}$.
Assume that
$$\limsup_{x\to 0}{{J(x)}\over {\|x\|^2}}=+\infty\ .$$
Then, for each $\rho>0$, there exists an open interval $I\subseteq ]0,+\infty[$ and an increasing function
$\varphi:I\to ]0,\rho[$ such that, for each $\lambda\in I$, one has
$$\emptyset\neq \left \{x\in S_{\varphi(\lambda)} : J(x)=\sup_{S_{\varphi(\lambda)}}J\right\}\subseteq
\{x\in X : x=\lambda J'(x)\}\ .$$
\bigskip
\bigskip
\bigskip
\bigskip
Here and in what follows,  $X$ is a real Hilbert space and $J:X\to {\bf R}$ is a functional, with $J(0)=0$.
For each $r>0$, set
$$S_{r}=\{x\in X : \|x\|^2=r\}\ ,$$
$$B_{r}=\{x\in X : \|x\|^2\leq r\}\ .$$
A point $\hat x\in S_{r}$ such that
$$J(\hat x)=\sup_{S_{r}}J$$
is called a spherical maximum of $J$. Assuming that $J$ is $C^1$, spherical maxima are important in
connection with the eigenvalue problem
$$J'(x)=\mu x\ .\eqno{(1)}$$
Actually, if $\hat x$ is a spherical maximum of $J$, by the classical Lagrange
multiplier theorem, there exists $\mu_{\hat x}\in {\bf R}$ such that
$$J'(\hat x)=\mu_{\hat x}\hat x\ .$$
More specifically, one could be interested in the multiplicity of solutions for $(1)$,
in the sense of finding some $\mu\in {\bf R}$ for which  there are more points
$x$ satisfying $(1)$. In this connection, however, just because of dependence
of $\mu_{\hat x}$ on $\hat x$, the existence of more spherical maxima in $S_{r}$
 does not imply automatically the existence of some $\mu\in {\bf R}$ for which
$(1)$ has more solutions. So, in order to the multiplicity of solutions of $(1)$, it is important
to know when, at least for some $r>0$,
 the spherical maxima in $S_{r}$ share the same Lagrange multiplier.\par
\smallskip
The aim of the present note is to give a contribution along such a direction. \par
\smallskip
 Here is our basic result:\par
\medskip
THEOREM 1.  - {\it  For some $\rho>0$, assume that
$J$ is  G\^ateaux differentiable in \hbox {\rm int}$(B_{\rho})\setminus \{0\}$ 
and that
$${{\beta_{\rho}}\over {\rho}}<\delta_{\rho} \eqno{(2)}$$
where
$$\beta_{\rho}=\sup_{B_{\rho}}J$$
and
$$\delta_{\rho}=\sup_{x\in B_{\rho}\setminus\{0\}}{{J(x)}\over
{\|x\|^2}}\ .$$
Assume also that, for some $a>0$, with
$$a>{{\rho}\over {\rho\delta_{\rho}-\beta_{\rho}}}$$
if $\delta_{\rho}<+\infty$, 
the restriction of the functional $\|\cdot\|^2-a J(\cdot)$ to $B_{\rho}$ is sequentially
weakly lower semicontinuous.\par
For each $r\in ]\beta_{\rho},+\infty[$, put
$$\eta(r)=\sup_{y\in B_{\rho}}{{\rho-\|y\|^2}\over {r-J(y)}}$$
and
$$\Gamma(r)=\left \{ x\in B_{\rho} : {{\rho-\|x\|^2}\over {r-J(x)}}=\eta(r)\right \}\ .$$
Then, the following assertions hold:\par
\noindent
$(i)$\hskip 5pt the function $\eta$ is convex and decreasing in $]\beta_{\rho},+\infty[$, with $\lim_{r\to +\infty}\eta(r)=0$\ ;\par
\noindent
$(ii)$\hskip 5pt for each $r\in \left ] \beta_{\rho}+{{\rho}\over {a}},\rho\delta_{\rho}\right [ $, 
the set $\Gamma(r)$ is non-empty and, for
every $\hat x\in \Gamma(r)$, one has
$$0<\|\hat x\|^2<\rho$$
and
$$\hat x\in \left \{ x\in S_{\|\hat x\|^2} : J(x)=\sup_{S_{\|\hat x\|^2}}J\right \} \subseteq
\left \{ x\in \hbox {\rm int}(B_{\rho}) : \|x\|^2-\eta(r)J(x)=\inf_{y\in B_{\rho}}(\|y\|^2-\eta(r)J(y))\right \}$$
$$\subseteq \left \{ x\in X : x={{\eta(r)}\over {2}} J'(x)\right \}\ ;$$
$(iii)$\hskip 5pt for each $r_1, r_2\in \left ] \beta_{\rho}+{{\rho}\over {a}},\rho\delta_{\rho}\right [$, with $r_1<r_2$, and
each $\hat x\in \Gamma(r_1), \hat y\in \Gamma(r_2)$, one has
$$\|\hat y\|<\|\hat x\|\ ;$$
$(iv)$\hskip 5pt if $A$ denotes
the set of all  $r\in \left ]\beta_{\rho}+{{\rho}\over
{a}},\rho\delta_{\rho}\right [$ such that $\Gamma(r)$ is a singleton, then
the function $r\to \Gamma(r)$ ($r\in A$) is continuous with respect to the weak topology; if, in
addition, $J$ is sequentially weakly upper semicontinuous in $B_{\rho}$, then $\Gamma_{|A}$
is continuous with respect to the strong topology.}
\smallskip
Before proving Theorem 1, let us recall a proposition from [1] that will be used in the proof:\par
\medskip
PROPOSITION 1. - {\it Let $Y$ be a non-empty set, 
$f, g:Y\to {\bf R}$ two functions, and $a, b$ two
real numbers, with $a<b$. Let $y_a$ 
 be a global minimum of the function $f+ag$ and
$y_b$ a global minimum of the function
$f+bg$.\par
Then, one has $g(y_b)\leq g(y_a)$.}\par
\medskip
{\it Proof of Theorem 1.} By definition, the function $\eta$ is the upper envelope of a family of functions
which are decreasing and convex in $]\beta_{\rho},+\infty][$. So, $\eta$ is convex and non-increasing. 
We also have
$$\eta(r)\leq {{\rho}\over {r-\beta_{\rho}}}\eqno{(3)}$$
for all $r>\beta_{\rho}$ and so
$$\lim_{r\to +\infty}\eta(r)=0\ .$$
In turn, this implies that $\eta$ is decreasing as it never vanishes. Now,
fix $r\in\left ]\beta_{\rho}+{{\rho}\over {a}}, \rho\delta_{\rho}\right [$. So, we have
$${{\rho}\over {r-\beta_{\rho}}}<a\ .$$
Consequently, by $(3)$, 
$$\eta(r)<a\ .$$
Observe that, for each $\lambda\in ]0,a[$, the restriction
to $B_{\rho}$ of the functional $\|\cdot\|^2-\lambda J(\cdot)$
is sequentially weakly lower semicontinuous. In this connection,
it is enough to notice that
$${{a}\over {a-\lambda}}(\|x\|^2-\lambda J(x))=
\|x\|^2+{{\lambda}\over {a-\lambda}}(\|x\|^2-a J(x))\ .$$
Fix a sequence $\{x_n\}$ in $B_{\rho}$ such that
$$\lim_{n\to \infty}{{\rho-\|x_n\|^2}\over {r-J(x_n)}}=\eta(r)\ .$$
Up to a sub-sequence, we can suppose that $\{x_n\}$ converges weakly
to some $\hat x_r\in B_{\rho}$.  Fix $\epsilon\in ]0,\eta(r)[$. 
For each $n\in {\bf N}$ large enough, we have
$${{\rho-\|x_n\|^2}\over {r-J(x_n)}}>\eta(r)-\epsilon$$
and so
$$\|x_n\|^2+(\eta(r)-\epsilon)(r-J(x_n))<\rho\ .$$
But then, by sequential weak lower semicontinuity, we have
$$\|\hat x_r\|^2+(\eta(r)-\epsilon)(r-J(\hat x_r))\leq \liminf_{n\to \infty}
(\|x_n\|^2+(\eta(r)-\epsilon)(r-J(x_n)))\leq\rho\ .$$
Hence, since $\epsilon$ is arbitrary, we have
$$\|\hat x_r\|^2+\eta(r)(r- J(\hat x_r))\leq \rho$$
and so
$${{\rho-\|\hat x_r\|^2}\over {r-J(\hat x_r)}}=\eta(r)\ ,$$
that is $\hat x_r\in \Gamma(r)$. Now, let $\hat x$ be any point of $\Gamma(r)$.
Let us show that $\hat x\neq 0$. Indeed, since ${{r}\over {\rho}}<\delta_{\rho}$,
 there exists
$\tilde x\in B_{\rho}\setminus \{0\}$ such that
$${{J(\tilde x)}\over {\|\tilde x\|^2}}>{{r}\over {\rho}}\ .$$
Clearly, this is equivalent to
$${{\rho}\over {r}}<{{\rho-\|\tilde x\|^2}\over {r-J(\tilde x)}}\ .$$
So
$${{\rho}\over {r}}<{{\rho-\|\hat x\|^2}\over {r-J(\hat x)}}$$
and hence, since $J(0)=0$, we have $\hat x\neq 0$, as claimed. Clearly,
$\|\hat x\|^2<\rho$ as $\eta(r)>0$. Moreover, if $x\in S_{\|\hat x\|^2}$,
we have
$${{1}\over {r-J(x)}}\leq {{1}\over {r-J(\hat x)}}$$
from which we get
$$J(\hat x)=\sup_{S_{\|\hat x\|^2}}J\ .$$
Now, let $u$ be any global maximum of $J_{|S_{\|\hat x\|^2}}$. Then, we have
$${{\rho-\|u\|^2}\over {r-J(u)}}=\eta(r)$$
and so
$$\|u\|^2-\eta(r) J(u)=\rho-r\eta(r)\leq \|x\|^2-\eta(r) J(x)$$
for all $x\in B_{\rho}$. Hence, as $\|u\|^2<\rho$, the point
$u$ is a local minimum of the functional $\|\cdot\|^2-\eta(r) J(\cdot)$.
Consequently, we have
$$u={{\eta(r)}\over {2}} J'(u)\ ,$$
and the proof of $(ii)$ is complete. To prove $(iii)$, observe that
$${{1}\over {\eta(r)}}=\inf_{\|x\|^2<\rho}{{r-J(x)}\over {\rho-\|x\|^2}}\ .$$
As a consequence, for each 
$r_1, r_2\in ]\beta_{\rho}+{{\rho}\over {a}},\rho\delta_{\rho}[$,
with $r_1<r_2$,  and for each $\hat x\in \Gamma(r_1), \hat y\in \Gamma(r_2)$, we have
$${{r_1-J(\hat x)}\over {\rho-\|\hat x\|^2}}=
\inf_{\|x\|^2<\rho}{{r_1-J(x)}\over {\rho-\|x\|^2}}$$
and
$${{r_2-J(\hat y)}\over {\rho-\|\hat y\|^2}}=
\inf_{\|x\|^2<\rho}{{r_2-J(x)}\over {\rho-\|x\|^2}}\ .$$
Therefore, in view of Proposition 1, we have
$${{1}\over {\rho-\|\hat y\|^2}}\leq {{1}\over {\rho-\|\hat x\|^2}}$$
and so
$$\|\hat y\|\leq \|\hat x\|\ .$$
We claim that 
$$\|\hat y\|<\|\hat x\|\ .$$
Arguing by contradiction, assume that $\|\hat y\|=\|\hat x\|$. In view of $(ii)$, this would imply that
$J(\hat y)=J(\hat x)$ and so, at the same time, 
$$\hat y={{\eta(r_2)}\over {2}}J'(\hat y)$$
and
$$\hat y={{\eta(r_1)}\over {2}}J'(\hat y)\ .$$
In turn, this would imply $\eta(r_1)=\eta(r_2)$ and hence $r_1=r_2$, a contradiction. So,
$(iii)$ holds. Finally, let us prove $(iv)$. 
For each $r\in A$,  continue to denote by $\Gamma(r)$
the unique point of $\Gamma(r)$.
Let $r\in A$ and let
$\{r_k\}$ be any sequence in $A$
converging to
$r$. Up to a subsequence, $\{\Gamma(r_k)\}$ converges
weakly to some $\tilde x\in B_{\rho}$. 
Moreover, for each $k\in {\bf N}$, $x\in B_{\rho}$, one has
$${{\rho-\|x\|^2}\over {r_k-J(x)}}\leq {{\rho-\|\Gamma(r_k)\|^2}\over {r_k-J(\Gamma(r_k))}}\ .$$
From this, after easy manipulations, we get
$$\|\Gamma(r_k)\|^2-{{\rho-\|x\|^2}\over {r-J(x)}}J(\Gamma(r_k))-
\left ( {{\rho-\|x\|^2}\over {r_k-J(x)}}-{{\rho-\|x\|^2}\over {r-J(x)}}\right )J(\Gamma(r_k))\leq
\rho-{{\rho-\|x\|^2}\over {r_k-J(x)}}r_k\ .\eqno{(4)}$$
Since the sequence $\{J(\Gamma(r_k)\}$ is bounded above, we have
$$\limsup_{k\to \infty}\left ( {{\rho-\|x\|^2}\over {r_k-J(x)}}-{{\rho-\|x\|^2}\over {r-J(x)}}\right )J(\Gamma(r_k))\leq 0\ .\eqno{(5)}$$
On the other hand, by sequential weak semicontinuity, we also have
$$\|\tilde x\|^2-{{\rho-\|x\|^2}\over {r-J(x)}}J(\tilde x)\leq
\liminf_{k\to \infty}\left ( \|\Gamma(r_k)\|^2-{{\rho-\|x\|^2}\over {r-J(x)}}J(\Gamma(r_k))\right )\ .\eqno{(6)}$$
Now, passing in $(4)$ to the $\liminf$, in view of $(5)$ and $(6)$, we obtain
$$\|\tilde x\|^2-{{\rho-\|x\|^2}\over {r-J(x)}}J(\tilde x)\leq
\rho-{{\rho-\|x\|^2}\over {r-J(x)}}r$$
which is equivalent to
$${{\rho-\|x\|^2}\over {r-J(x)}}\leq  {{\rho-\|\tilde x\|^2}\over
{r-J(\tilde x)}}\ .$$
Since this holds for all $x\in B_{\rho}$, we have $\tilde x=\Gamma(r)$. So, $\Gamma_{|A}$ is continuous
at $r$ with respect to the weak topology. Now, assuming also that $J$ is  sequentially weakly upper semicontinuous,
in view of the continuity of $\eta$ in $]\beta_{\rho},+\infty[$, we have
$$\lim_{k\to \infty} {{\rho-\|\Gamma(r_k)\|^2}\over {r_k-J(\Gamma(r_k))}}=
{{\rho-\|\Gamma(r)\|^2}\over {r-J(\Gamma(r))}}\ ,$$
and hence
$$\liminf_{k\to \infty}(\rho-\|\Gamma(r_k)\|^2)=
{{\rho-\|\Gamma(r)\|^2}\over {r-J(\Gamma(r))}}
\liminf_{k\to \infty}(r_k-J(\Gamma(r_k)))=
{{\rho-\|\Gamma(r)\|^2}\over {r-J(\Gamma(r))}}
(r-\limsup_{k\to \infty}J(\Gamma(r_k)))$$
$$\geq {{\rho-\|\Gamma(r)\|^2}\over {r-J(\Gamma(r))}}
(r-J(\Gamma(r))=
\rho-\|\Gamma(r)\|^2$$
from which
$$\limsup_{k\to \infty}\|\Gamma(r_k)\|\leq \|\Gamma(r)\|\ .$$
 Since $X$ is a Hilbert space and $\{\Gamma(r_k)\}$ converges weakly to $\Gamma(r)$, this implies that
$$\lim_{k\to \infty}\|\Gamma(r_k)-\Gamma(r)\|=0\ ,$$
which shows the continuity of $\Gamma_{|A}$ at $r$ in the strong topology.\hfill $\bigtriangleup$
\medskip
REMARK 1. - Clearly, when $J$ is sequentially weakly upper semicontinuous in $B_{\rho}$, 
the assertions of Theorem 1 hold in the whole interval $]\beta_{\rho},\rho\delta_{\rho}[$, since
$a$ can be any positive number.\par
\medskip
REMARK 2. - The simplest way to satisfy condition $(2)$ is, of course, to assume that
$$\limsup_{x\to 0}{{J(x)}\over {\|x\|^2}}=+\infty\ .$$
Another reasonable way is provided by the following proposition:\par
\medskip
PROPOSITION 2. - {\it For some $s>0$, assume that $J$ is G\^ateaux differentiable in
$B_s\setminus \{0\}$ and that there exists a global maximum $\hat x$
of $J_{|B_s}$ such that
$$\langle J'(\hat x),\hat x\rangle <2J(\hat x)\ .$$
Then, $(2)$ holds with $\rho=\|\hat x\|^2$\ .}\par
\smallskip
PROOF. For each $t\in ]0,1]$, set
$$\omega(t)={{J(t\hat x)}\over {\|t\hat x\|^2}}\ .$$
Clearly, $\omega$ is derivable in $]0,1]$. In particular, one has
$$\omega'(1)={{\langle J'(\hat x),\hat x\rangle - 2J(\hat x)}\over {\|\hat x\|^2}}\ .$$
So, by assumption, $\omega'(1)<0$ and hence, in a left neighbourhood of $1$, we have
$$\omega(t)>\omega(1)$$
which implies the validity of $(2)$ with $\rho=\|\hat x\|^2$.\hfill $\bigtriangleup$\par
\medskip
Also, notice the following consequence of Theorem 1:\par
\medskip
THEOREM 2. - {\it For some $\rho>0$, let the assumptions of Theorem 1 be satisfied.\par
Then, there exists an open interval $I\subseteq ]0,+\infty[$ and an increasing function
$\varphi:I\to ]0,\rho[$ such that, for each $\lambda\in I$, one has
$$\emptyset\neq \left \{x\in S_{\varphi(\lambda)} : J(x)=\sup_{S_{\varphi(\lambda)}}J\right\}\subseteq
\{x\in X : x=\lambda J'(x)\}\ .$$}
PROOF. Take
$$I={{1}\over {2}}\eta\left (\left ]\beta_{\rho}+{{\rho}\over {a}},\rho\delta_{\rho}\right [\right )\ .$$
Clearly, $I$ is an open interval since $\eta$ is continuous and decreasing. Now, for each
$r\in \left ]\beta_{\rho}+{{\rho}\over {a}},\rho\delta_{\rho}\right [$, pick $v_r\in \Gamma(r)$. Finally,
set
$$\varphi(\lambda)=\|v_{\eta^{-1}(2\lambda)}\|^2$$
for all $\lambda\in I$. Taking $(iii)$ into account, we then realize that the function
$\varphi$ (whose range is contained in $]0,\rho[$)
is the composition of two decreasing functions, and so it is increasing. Clearly,
the conclusion follows directly from $(ii)$.\hfill $\bigtriangleup$
\medskip
We conclude deriving from Theorem 1 the following multiplicity result:\par
\medskip
THEOREM 3. - {\it For some $\rho>0$, assume that $J$ is sequentially weakly upper
semicontinuous in $B_{\rho}$ , G\^ateaux differentiable in $\hbox {\rm int}(B_{\rho})\setminus \{0\}$
and satisfies $(2)$. Moreover, assume that there exists $\tilde\rho$
satisfying
$$\inf_{x\in D}\|x\|^2<\tilde\rho<\sup_{x\in D}\|x\|^2\ ,\eqno{(7)}$$
where
$$D=\bigcup_{r\in ]\beta_{\rho},\rho\delta_{\rho}[}\Gamma(r)\ ,$$
such that $J_{|S_{\tilde\rho}}$ has either two global maxima or a global maximum at which
$J'$ vanishes.\par
Then, there exists $\tilde\lambda>0$ such that the equation
$$x=\tilde\lambda J'(x)$$
has at least two non-zero solutions  which are global minima of the restriction of the functional
 ${{1}\over {2}}\|\cdot\|^2-\tilde\lambda  J(\cdot)$ to $\hbox {\rm int}(B_{\rho})$.}\par
\smallskip
PROOF. For each $r\in ]\beta_{\rho},\rho\delta_{\rho}[$, in view of $(7)$,
we can pick  $v_r\in \Gamma(r)$ (recall
Remark 1), so that
$$\inf_{]\beta_{\rho},\rho\delta_{\rho}[}\psi<\tilde\rho<\sup_{]\beta_{\rho},\rho\delta_{\rho}[}\psi\ ,\eqno{(8)}$$
where
$$\psi(r)=\|v_r\|^2\ .$$
Two cases can occur. First, assume that $\tilde\rho\in \psi( ]\beta_{\rho},\rho\delta_{\rho}[)$.
So, $\psi(r)=\tilde\rho$ for some $\tilde r\in  ]\beta_{\rho},\rho\delta_{\rho}[$. So, by (ii), for each
global maximum $u$ of $J_{|S_{\tilde\rho}}$, we have $J'(u)\neq 0$.
 As a consequence, in this case,
$J_{|S_{\tilde\rho}}$ has at least two global maxima which, by $(ii)$ again, satisfies the conclusion
with $\tilde\lambda={{1}\over {2}}\eta(\tilde r)$. Now, suppose that 
$\tilde\rho\not\in \psi( ]\beta_{\rho},\rho\delta_{\rho}[)$. In this case, in view of $(8)$, 
the function $\psi$ is discontinuous and hence, in view of $(iv)$, there exists some
$r^*\in ]\beta_{\rho},\rho\delta_{\rho}[$ such that $\Gamma(r^*)$ has at least two elements which,
by $(ii)$, satisfy the conclusion with $\tilde\lambda={{1}\over {2}}\eta(r^*)$.\hfill $\bigtriangleup$\par
\bigskip
\centerline {\bf References}\par
\bigskip
\noindent
[1]\hskip 5pt B. RICCERI, {\it Uniqueness properties of functionals with
Lipschitzian derivative}, Port. Math. (N.S.), {\bf 63} (2006), 393-400.\par
\bigskip
\bigskip
Department of Mathematics\par
University of Catania\par
Viale A. Doria 6\par
95125 Catania\par
Italy\par
{\it e-mail address}: ricceri@dmi.unict.it

\bye